\documentstyle[11pt]{article}

\newenvironment{beweis}{{\it Proof.}\ }{\ $\ \ \ \Diamond$ \\\ }

 \newcounter{nsatz}[section]
 \newcounter{nlemma}[section]
 \newcounter{ndef}[section]
 \newcounter{nhyp}[section]
 \newcounter{nconjecture}[section]
 \newcounter{ncor}[section]
 \newcounter{nrem}[section]
 \newcounter{nexample}[section]
 \newcounter{nprop}[section]

 \newenvironment{nsatz}{\refstepcounter{nsatz}{\bf \arabic{section}.\arabic{nsatz}}\
                {\sc\bf Theorem.\ }\it}{\\\\ \rm}

 \newenvironment{nlemma}{\setcounter{nlemma}{\value{nsatz}}
                \refstepcounter{nlemma}
                \setcounter{nsatz}{\value{nlemma}}
                {\bf \arabic{section}.\arabic{nsatz}}\
                {\sc\bf Lemma.\ }\it}{\\\\ \rm}

 \newenvironment{ncor}{\setcounter{ncor}{\value{nsatz}}
                \refstepcounter{ncor}
                \setcounter{nsatz}{\value{ncor}}
                {\bf \arabic{section}.\arabic{nsatz}}\
                {\sc\bf Corollary.\ }\it}{\\\\ \rm}

\begin{document}
\newcommand{\n}{{\mbox{\rm I$\!$N}}}
\newcommand{\z}{{\mbox{{\sf Z\hspace{-0.4em}Z}}}}
\newcommand{\R}{{\mbox{\rm I$\!$R}}}
\newcommand{\Q}{{\mbox{\rm I$\!\!\!$Q}}}
\newcommand{\C}{{\mbox{\rm I$\!\!\!$C}}}
\newcommand{\ug}{\ \raisebox{-.3em}{$\stackrel{\scriptstyle \leq}
{\scriptstyle \sim}$} \ }
\newcommand{\X}{{\mbox{$\setminus$}\mbox{$\!\!\!/$}}}
\thispagestyle{empty}
\setlength{\parindent}{0pt}
\setlength{\parskip}{5pt plus 2pt minus 1pt}

\thispagestyle{empty}
\newcommand{\Syl}{{\mbox{\rm Syl}}}
\newcommand{\Hall}{{\mbox{\rm Hall}}}
\newcommand{\cl}{{\mbox{\rm cl}}}
\newcommand{\Irr}{{\mbox{\rm Irr}}}
\newcommand{\GL}{{\mbox{\rm GL}}}
\newcommand{\SL}{{\mbox{\rm SL}}}
\newcommand{\Ortho}{{\mbox{\rm O}}}

\newcommand{\sIrr}{{\mbox{\scriptsize\rm Irr}}}
\newcommand{\Char}{{\mbox{\rm Char}}}
\mbox{\vspace{4cm}}
\vspace{3cm}
\begin{center}
{\bf \Large\bf  Class 2 quotients of solvable linear groups\\}
\vspace{2cm}
by\\
\vspace{11pt}
Thomas Michael Keller\\
Department of Mathematics\\
Texas State University\\
601 University Drive\\
San Marcos, TX 78666\\
USA\\
E--mail: keller@txstate.edu\\
\vspace{1cm}
and\\
\vspace{1cm}

Yong Yang\\
Department of Mathematics\\
Texas State University\\
601 University Drive\\
San Marcos, TX 78666\\
USA\\
E--mail: yang@txstate.edu\\
\vspace{1cm}
\vspace{1cm}
\end{center}
\thispagestyle{empty}
\newpage

\begin{center}
\parbox{12.5cm}{{\small
{\sc Abstract.}
Let $G$ be a finite group, and let $V$ be a completely reducible faithful $G$-module. By a result of Glauberman it has been
known for a long time that if $G$ is nilpotent of class 2, then $|G| < |V|$. In this paper we
generalize this result as follows. Assuming $G$ to be solvable, we show that the order of the
maximal class 2 quotient of $G$
is strictly bounded above by $|V|$.}}
\end{center}
\normalsize

\section{Introduction}\label{section0}

In \cite[Proposition 1]{glauberman} G. Glauberman proved that if $n$ is a positive integer,
$p$ is a prime and $G\leq \mbox{GL}(n,p)$
is a $p'$-group which is nilpotent of class 2, then $|G| < |V|$.

The goal of this paper is to generalize this result as follows. For a finite group $G$
put $G^c=[G,G,G]=[G', G]$; i.e., $G^c$ is the intersection of all normal
subgroups of $G$ whose quotient group is nilpotent of class 2, so that $G/G^c$ is the (maximal) class 2 quotient
of $G$. With this we will prove the following generalization of Glauberman's result.\\

\begin{nsatz}\label{mainthm}
Let $G$ be a finite solvable group and $V\ne 0$ a finite faithful completely reducible $G$-module,
possibly of mixed characteristic. Then
\[|G/G^c|\ <\  |V|.\]
\end{nsatz}
Of course, when $G$ is nilpotent of class 2, this is just Glauberman's result, which we will use in the proof of the
above theorem.\\

\ref{mainthm} can also be viewed as a strengthening - for solvable groups and completely reducible modules - of
a result by Aschbacher and Guralnick, see \cite[Theorem 1]{aschbacher-guralnick}. They proved that the order of the
abelian quotient $|G/G'|$ of $G$, i.e., the class 1 quotient of $G$, is strictly bounded above by $|V|$, where
$G$ is a finite faithful linear group on the finite module $V$ such that O$_r(G)=1$ for the characteristic $r$ of $V$.
For solvable $G$ and completely reducible $V$
our new result shows that even the class 2 quotient of $G$ is strictly bounded above by $|V|$.\\

We note that we believe that the main result of this paper remains true for arbitrary finite groups in place of solvable groups. \\

\section{Bounding the class 2 quotient by the module size}\label{section1}

In this section we prove \ref{thm1}. We first prove a reduction lemma.\\

\begin{nlemma}\label{lem1}
Let $G$ be a finite group and $N\unlhd G$. Then
\[|G/G^c|=|G/G^cN|\cdot|N:N\cap G^c|;  \]
and
\[|G:G^c|\mbox{ divides }|G/N:(G/N)^c||N:N^c|.\]
\end{nlemma}
\begin{beweis}
The proof goes along exactly the same lines as the proof of \cite[Lemma1]{keller-yang}
\end{beweis}

In \cite{keller-yang} the following result is proved, albeit not stated as a separate result.\\

\begin{nlemma}\label{lem2}
Let $G$ be a finite solvable group and $V$ a finite faithful completely reducible $G$-module,
possibly of mixed characteristic. If $G$ is not nilpotent, then
\[|G/G'|\ <\  |F(G)/F(G)'|.\]
\end{nlemma}
\begin{beweis}
This is not explicitly stated in \cite{keller-yang}, but is contained in the proof of \cite[Theorem 2.3]{keller-yang}
at the beginning where it is obtained that $G$ can be assumed to be nilpotent. Note that in that proof the
inequality (2) there is obtained by using an inductive hypothesis which we do not have if we just want to prove the
current lemma. But we can, instead, simply use the now proven \cite[Theorem 2.3]{keller-yang} to get that inequality (2)
in the proof of the current lemma, or, alternatively, use the Aschbacher-Guralnick bound \cite[Theorem 3]{aschbacher-guralnick}
to obtain (2), and otherwise just follow the argument outlined in the proof of \cite[Theorem 2.3]{keller-yang}.
In particular, if one uses the Aschbacher-Guralnick result, then the proof is independent of \cite[Theorem 2.3]{keller-yang}.
\end{beweis}

We sort of used \ref{lem2} in \cite{keller-yang} to reduce the proof of the main result to nilpotent groups. We plan to do
a similar approach in the proof of the main result \ref{thm1} of this paper, and as a consequence, will get the class 2-analogue
of \ref{lem2} (see \ref{cor1} below); unlike in the proof of \ref{lem2}, however, we do not see how to replace using \ref{thm1} by some other result,
so that the class 2-analogue of \ref{lem2} truly appears to be a consequence \ref{thm1}.\\

We therefore now start proving \ref{mainthm} and roughly use a similar structure as in \cite[Theorem 2.3]{keller-yang}.
Not surprisingly, the reduction to nilpotent groups is a little bit more involved, but the remaining proof is a little shorter, since we can rely
on work of Glauberman \cite{glauberman}.\\

\begin{nsatz}\label{thm1}
Let $G$ be a finite solvable group and $V \ne 0$ a finite faithful completely reducible $G$-module, possibly of mixed characteristic.
Then
\[|G/G^c|\ <\ |V|.\]
\end{nsatz}
\begin{beweis}
We work by induction on $|GV|$.
First we want to reduce to the case that $G$ is nilpotent.\\

So suppose $G$ is not nilpotent. Write $F=F(G)$ for the Fitting subgroup and $\Phi=\Phi(G)$ for
the Frattini subgroup of $G$. As $G$ is not nilpotent, we have $F<G$. \\
Moreover, by Gasch\"{u}tz' theorem $F/\Phi$ is a faithful, completely reducible $G/F$-module
(possibly of mixed characteristic). We now write $F/\Phi=W_1\oplus W_2$, where
$W_1=(F\cap G^c)\Phi/\Phi$ and $W_2$ is a $G$-invariant complement of $W_1$ in $F/\Phi$. Hence
\[W_2\cong(F/\Phi)/W_1\cong F/(F\cap G^c\Phi)\mbox{ as
}G\mbox{-modules}.\]

We now claim that $G/F$ acts faithfully on $W_1$. Assume not, then there exists a $g\in G-F$ such that
$g$ is of prime order $q$ for some prime $q$ and $g$ acts trivially on $W_1$. Since by Gasch\"{u}tz' theorem
we know that there exists a subgroup $\Phi\leq L\leq G$ such that $L/\Phi$ is a complement of $F/\Phi$ in $G/\Phi$,
by choosing $g$ in $L$ we may assume that $g^q\in\Phi$.

 Since $g$ acts nontrivially
on $F/\Phi$, it follows that $g$ acts nontrivially on $W_2$. Then there exists an irreducible $G/F$-submodule $X_2$ of $W_2$
such that $g$ acts nontrivially on $X_2$. Then $X_2$ is of characteristic $p$ for a prime $p$, and there exists
an $x\in F$ of $p$-power order such that $x\Phi\in X_2$ and $[g,x]\in F-\Phi$. Now if $[g,x]\in G^c\Phi$, then $[g,x]\in F\cap G^c\Phi=(F\cap G^c)\Phi$
(the latter equality following from Dedekind's identity) and hence
$1\ne [g,x\Phi]\in X_2 \cap (F\cap G^c)\Phi/\Phi\leq W_2\cap W_1=1$, a contradiction. This shows that
$[x,g]\not\in G^c\Phi\ $ (+).
In particular, $x$, $g\in G-G^c\Phi$.\\

We consider two cases.\\

Case 1: It is possible to choose $X_2$ in such a way that $p\ne q$.\\
Then we do this, so $p\ne q$. Now since $G/G^c$ is nilpotent, $p\ne q$, and $xG^c\in G/G^c$ is of $p$-power order, we see that
$1=[xG^c, gG^c]=[x,g]G^c$ which implies that $[x,g]\in G^c$, contradicting (+). Hence this case cannot occur.\\

Case 2: It is not possible to choose $X_2$ in such a way that $p\ne q$.\\
Then $p=q$ and $g$ acts trivially on the Hall $q'$-subgroup of $F$. (More precisely, in Case 2 first we know that
$g$ acts trivially on the Hall $q'$-subgroup of $F/\Phi$, but well-known results on coprime automorphisms then imply
that $g$ indeed acts trivially on the Hall $q'$-subgroup of $F$.)\\

Recall that $g\not\in F$, and we now work towards a contradiction. This contradiction will show that Case 2 cannot occur either
and thus we will have shown that $G/F$ acts faithfully on $W_1$.\\

Now let $H$ be a Hall $p'$-subgroup of $F$ and put $C=C_G(H)$. Clearly $H\unlhd G$ and $C\unlhd G$. Also, since $F\leq C$,
it follows that $F(C)=F$. Hence $g\in C-F(C)$. We also have $O_p(C/F)=1$, because if $P_0$ is the inverse image of
$O_p(C/F)$ in $C$ and $P$ is a Sylow $p$-subgroup of $C$, then we have $P_0=P\times H$, so $P_0$ is nilpotent and normal in $C$
and hence $P_0\leq F(C)=F$ and so $O_p(C/F)=1$ as claimed. As $g$ is of order $p$, we obtain that $g\not\in F_2(G)$ (the
second Fitting subgroup of $G$).\\
Hence there exists a prime $r\ne p$ such that $g$ acts nontrivially on the Sylow $r$-subgroup $R_0$ of $F(C/F)$. By the
Hall-Higman reduction \cite[III, Satz 13.5]{huppert} we can find a subgroup $R_1$ of $R_0$ which is minimal with respect to $g$
acting nontrivially on it, and $R_1$ is special and $[\langle g\rangle, R_1]=R_1$.\\
We observe that since $G$ acts trivially on $W_1$, we have $[\langle g\rangle, W_1]=1$ and thus
$[W_1, \langle g\rangle, R_1]=1$ and also $[R_1, W_1,\langle g\rangle]=1$. So by the Three-subgroups-lemma we conclude that
$[R_1, W_1]=[\langle g\rangle, R_1, W_1]=1$, that is, $R_1$ acts trivially on $W_1$ and thus nontrivially on $W_2$. \\

Since $[\langle g\rangle, R_1]=R_1$, it follows that $g$ will act nontrivially on any submodule of $W_2$ on which $R_1$
acts nontrivially. In particular, by possibly replacing $X_2$ we may assume that $R_1$ acts nontrivially on $X_2$.\\

Now let $R\leq G$ be the inverse image of $R_1\leq C/F$ in $C$, and put $T=\langle g \rangle R$. Moreover let $Y\leq X_2$ be an irreducible $T$-module on which $R_1$ acts nontrivially. As seen above, then $g$ will act
nontrivially on $Y$. Put $S=TY\leq G/\Phi$. We consider again two cases.\\

Case A: $Y\not\leq S^c$. \\
Since $Y$ is irreducible as $T$-module, this means that $S^c\cap Y=1$.
So if $R\leq S^c$, then $[R,Y]\leq S^c\cap Y=1$ contradicting $R$ acting nontrivially on $Y$.
Hence $R\not\leq S^c$, and clearly $g\Phi\not\in S^c$ (as $g\not\in G^c\Phi$ and $S^c\leq (G/\Phi)^c=G^c\Phi/\Phi$).
But by the choice of $R_1$ it then is clear that $S/S^c$ cannot be nilpotent since $g$ cannot act trivially on any
factor group of $R_1$). This contradiction shows that Case A cannot occur.\\

Case B: $Y\leq S^c$.\\
Let $y\in Y$ such that $[g,y]\ne 1$. Then $1\ne [g,y]\in Y$, and as before we conclude that
$[g,y]\in (F\cap G^c)\Phi/\Phi$ (since $S^c\leq G^c$) and thus $[g, y]\in \left((F\cap G^c)\Phi/\Phi\right)\cap X_2
\leq W_1\cap W_2=1$, a contradiction. So Case B is impossible, too.\\

So this finally refutes our assumption that $G/F$ does not act faithfully on $W_1$, and thus we now suppose that
$W_1$ is a faithful $G/F$-module.\\
Now by induction, applied to the action of $G/F$ on $W_1=(F\cap G^c)\Phi/\Phi$, we can conclude that
\[(1)\ \ \ |G/(G^cF)|\ =\ |(G/F)/(G^cF/F)|\ = \ |(G/F)/(G/F)^c|\  <\ |W_1|.\]
Furthermore
\[(2)\ \ \ |G/G^c|\ =\ |G/(G^cF)| \cdot |G^cF/G^c|\ =\ |G/(G^cF)| \cdot |F/(F\cap G^c)|\ <\
|W_1| \cdot |F/(F\cap G^c)|,\]
where the inequality follows from (1). Now as $F^c\leq \Phi(F)\leq \Phi(G)$, we see that
\[(3)\ \ \ |W_1|\ =\ |(F\cap G^c)\Phi/\Phi|\ =\ |(F\cap G^c)/(F\cap G^c\cap \Phi)|\ \leq\ |(F\cap G^c)/F^c|,\]
and this together with (2) yields
\[(4)\ \ \ |G/G^c|\ <\  |(F\cap G^c)/F^c| \cdot |F/(F\cap G^c)|\ =\ |F/F^c|.\]

Now since $F<G$, by induction we have $|F/F^c|< |V|$, and together with (4) we obtain $|G/G^c|< |V|$,
as desired.\\

So from now on we may assume that $G$ is nilpotent.\\

Next we show that we may assume that $V$ is irreducible. Suppose that
$V=V_1\oplus V_2$ where $V_1,V_2$ are nontrivial $G$-submodules of $V$. By induction, we have
that $|G/C_G(V_1):(G/C_G(V_1))^c| < |V_1|$. Since $C_G(V_1)$ is normal in $G$,
we know that $V$ is a completely reducible $C_G(V_1)$-module, and so $V_2$ is a faithful completely
reducible $C_G(V_1)$-module. Thus by induction, $|C_G(V_1):C_G(V_1)^c| < |V_2|$. So altogether with
Lemma 2.1 we obtain that
\[|G:G^c| \leq  |G/C_G(V_1):(G/C_G(V_1))^c|\ |C_G(V_1):C_G(V_1)^c| < |V_1||V_2|=|V|,\]
and we are done. So now suppose that $V$ is irreducible as $G$-module.\\

If $V$ is quasiprimitive, then it is well-known (see e.g. the proof of
\cite[Theorem 3.3]{manz-wolf}) that $G=S\times T$ where $T$ is cyclic of odd order and $S$ is a 2-group that is cyclic, dihedral,
quaternion or semi-dihedral. In particular, $G$ has a cyclic normal subgroup $U$ of index at most 2. Thus clearly $U$ has a regular
orbit on $V$ and so $|G|/2\ <\ |V|$. If $G$ is abelian, then even $|G|<|V|$, and we are done. If $G$ is not of class 2, then
$|G:G^c|\leq |G|/2\ <\ |V|$, and we are done. So we may assume that $G$ is of class 2. In this case $G$ must be of order 8 and
dihedral or quaternion, and in both cases it is easy to see that $|G|<|V|$ follows. Therefore from now on we may assume that
$V$ is not quasiprimitive.\\

Next we would like to further reduce to the case that $G$ is a $p$-group for some prime $p$. For this, we argue somewhat similarly as in
the corresponding part of the proof of \cite[Theorem 2.3]{keller-yang}, but have to work a little harder. As in \cite{keller-yang}
we will avoid using the main result in \cite{hargraves} to keep the argument here more elementary and self-contained; also, using
\cite{hargraves} would save us only little work.\\

First, working towards a contradiction, assume that there are at least
two distinct primes dividing $|G|$ for which $G$ has a nonabelian Sylow subgroup. Let $p$ be one of these primes and
write $G=P\times H$ where
$P\in\Syl_p(G)$ and $H\in\Hall_{p'}(G)$. Then both $P$ and $H$ are nonabelian. Now $V$ is a finite $G$-module over a finite field, so let $K$ be that field.
By \cite[Lemma 10]{robinson-thompson} there exists a finite field extension $L$ of $K$ such that if $U$ is any irreducible summand of
$V$ viewed a an $LG$-module, then the permutation actions of $G$ on $V$ and $U$ are permutation isomorphic. So by studying the
action of $G$ on $U$ instead of $V$, we may as well assume that $V$ is absolutely irreducible. Then by \cite[(3.16)]{aschbacher}
we may further assume that $V=X_1\otimes X_2$ is a tensor product, where $X_1$ is a faithful, irreducible $P$-module and $X_2$ is a faithful
irreducible $H$-module. So if $k_i=\dim_K X_i$ for $i=1,2$, then $k_i\geq 2$ and $|X_i|=|K|^k_i$ for $i=1,2$.
Now by induction, we know that $|P:P^c|< |X_1|$ and $|H:H^c|< |X_2|$, and therefore we get
\[|G:G^c|\ =\ |P:P^c|\cdot |H:H^c|\ <\ |X_1|\cdot |X_2|\ =\ |K|^{k_1+k_2}\ \leq |K|^{k_1k_2}\ =\ |V|,\]
and we are done. \\

Since $G$ is not abelian, we therefore know that there exists exactly one prime $p$ such that $G$ has a nonabelian Sylow $p$-subgroup.
Since $V$ is irreducible, every abelian Sylow subgroup of $G$ must be cyclic and act frobeniusly on $V$, and so we have  $G=P\times H$ where
$P\in\Syl_p(G)$ and $H\in\Hall_{p'}(G)$, and $H$ is cyclic and acts frobeniusly on $V$. To complete the reduction to $G$ being a $p$-group,
we eventually will show that
we may assume that $H=1$. \\

Since $V$ is not quasiprimitive, and $P$ is the only nonabelian Sylow subgroup of $G$, by \cite[Proposition 0.3]{manz-wolf}
there exists an
$E\unlhd G$ with $|G:E|=p$ such that $V_{E}=V_1\oplus\ldots\oplus V_p$ for homogeneous components $V_i$ of $V_{E}$
which are permuted in a $p$-cycle by $G/E$.\\
Clearly $H\leq E$, so write $E= D\times H$, where $D\in\Syl_p(E)$ is of index $p$ in $P$. \\

Now clearly $E^c=D^c$, and by the inductive hypothesis, $|E/D^c| < |V|$.\\
Next we want to show that we may assume that $D^c=G^c$. \\

To do this, assume that $D^c<G^c$.
Since $G$ is nilpotent, it is then clear that $p|D^c|=p|E^c|\leq |G^c|$, and so it follows that
\[|V|>|E/D^c|=\frac{|E|}{|D^c|}=\frac{p|E|}{p|D^c|}=\frac{|G|}{p|D^c|}\geq|G/G^c|,\]
and we are done. Thus from now on we assume that $D^c=G^c$.\\

Now put $C_0=D$, and if $C_i$ is already defined for some $i\in \{0,\dots ,p-1\}$, then define
\[C_{i+1}=C_{C_i}(V_{i+1})=C_i\cap C_{D}(V_{i+1})=\bigcap_{j=1}^{i+1}C_{E}(V_j).\]
Also let $D_i=C_i/C_{i+1}$
for $i\in\{0,\dots ,p-1\}$.
Since $C_{i+1}\unlhd C_i$, we see that
$V_{i+1}$ is a faithful completely reducible $D_i$-module for $i=0,\dots ,p-1$. Also note that
$D_{p-1}=\bigcap_{i=1}^{p-1}C_D(V_i)$ (which acts faithfully on $V_p$), and also $\prod\limits_{i=0}^{p-1}|D_i|=|D|$.\\

Now by induction we conclude that $|D_{i}:D_{i}^c|\ <\ |V_{i+1}|\ =\ |V_1|$ for $i=0,\dots ,p-1$.\\

We are now ready to show that $H=1$. For this, note that even $HD_i = (HC_i)/C_{i+1}$ acts faithfully and completely
reducibly on $V_{i+1}$, and so by induction for $E_i=(HC_i)/C_{i+1}$ we even get $|E_{i}:E_{i}^c|\ <\ |V_{i+1}|\ =\ |V_1|$ for $i=0,\dots ,p-1$.
In particular, clearly $|H|$ divides $|E_{i}:E_{i}^c|$ (since $E^c=D^c$), and hence we obtain $|D_{i}:D_{i}^c|\ <\ |V_{i+1}|/|H|\ =\ |V_1|$ for $i=0,\dots ,p-1$.
Consequently, we obtain
\[|G:G^c|=p|E:D^c|=p|H|\ |D:D^c|\leq p\ |H|\ \prod\limits_{i=0}^{p-1}|D_i:D_i^c|\ <\ p\cdot |H|\cdot \frac{|V|}{|H|^p}\ =\ \frac{p}{|H|^{p-1}}\ |V|.\]
So if $H>1$, then $|H|\geq 2$, and so $|G:G^c|< \frac{p}{2^{p-1}}|V|\leq |V|$ follows, and we are done.\\

Therefore for the rest of the proof we may assume that $G$ is a $p$-group.\\

Recall from above that $|D_{i}:D_{i}^c|\ <\ |V_{i+1}|\ =\ |V_1|$ for $i=0,\dots ,p-1$.\\

Next observe that there exists exactly one nonnegative integer $k$ such that
$|V_1|/p\leq p^k <|V_1|$.
Now if there exists a $j\in\{0,\dots ,p-1\}$ such that $|D_{j}:D_{j}^c|\ <\ p^k$, then $|D_{j}:D_{j}^c|\ <\ |V_1|/p$,
and then we obtain
\[|G:G^c|=p|D:D^c|\leq p\ \prod\limits_{i=0}^{p-1}|D_i:D_i^c|\ <\ p\cdot \frac{|V|}{p}\ =\ |V|,\]
and we are done.\\
So we may, therefore, assume from now on that $|D_{i}:D_{i}^c|\ =p^k$ for $i=0,\dots ,p-1$.\\

We next put $s=p^{1/p}$ and observe that if $p^k<|V_1|/s$, then
\[|G:G^c|=p|D:D^c|\leq p\ \prod\limits_{i=0}^{p-1}|D_i:D_i^c|\ <\ p\cdot \frac{|V|}{s^p}\ =\ p\cdot \frac{|V|}{p}\ =\ |V|,\]
and again we are done. Hence we now know that
\[|V_1|/s < p^k<|V_1|;\]
the first strict inequality follows from the fact that $|V_1|/s$ clearly is not an integer.\\

Now let \[K_i=\bigcap_{\j\in\{1,\dots ,p\}-\{i\}}C_D(V_j)\] ($i=1,\dots ,p$) and observe that $K_i\unlhd D$ for all $i$, and
$G/D$ permutes the $K_i$ transitively. Moreover, $K_p=D_{p-1}$, and thus $|K_i:K_i^c|=|D_{p-1}:D_{p-1}^c|=p^k>|V_1|/s$
for all $i$. Also, it is easy to see that $K:=\prod_{i=1}^{p}K_i=\X_{i=1}^p K_i$ is a direct product of the $K_i$ and
$K\unlhd G$ and thus acts completely reducibly on $V$. We therefore have
\[|K:K^c|\ =\ \prod\limits_{i=1}^{p}|K_i:K_i^c|\ >\ \prod\limits_{i=1}^{p}\frac{|V_1|}{s}\ =\ \frac{|V|}{p}.\]

Our next goal is to show that $D=K$. For this we first need to make another little observation.\\
Put $W_1=\sum\limits_{i=1}^{p-1}V_i$ and $W_2=V_p$. Then $D/K_p$ acts faithfully on $W_1$, and $K_p$ acts
faithfully on $W_2$, and both actions are clearly completely reducible. Now by \ref{lem1}
we have $|D/D^c|=|D/D^cK_p|\cdot|K_p:K_p\cap D^c|$, and by induction we have
$|D/D^cK_p|=|(D/K^p):(D/K_p)^c|<|W_1|$ and $|K_p:K_p\cap D^c|\leq |K_p:K_p^c|<|W_2|$. Now if even
$|K_p:K_p\cap D^c| < |K_p:K_p^c|$, then $p|K_p:K_p\cap D^c| \leq |K_p:K_p^c|$, and then
\begin{eqnarray*}
|G:G^c|=|G:D^c|=p|D:D^c|=p\cdot |D/D^cK_p|\cdot|K_p:K_p\cap D^c|\\
\leq |D/D^cK_p|\cdot|K_p:K_p^c|<|W_1|\cdot |W_2|=|V|,
\end{eqnarray*}
and we are done. Hence we may assume that $|K_p:K_p\cap D^c| = |K_p:K_p^c|$, or, equivalently, that
$K_p^c=K_p\cap D^c$. Since the $K_i$ are all conjugate in $G$, it follows that $K_i^c=K_i\cap D^c$ $(*)$ for
$i=1,\dots ,p$.\\

Now assume that $D>K$ and let $L\unlhd D$ such that $K<L$ and $|L:K|=p$.
Then $V$ is a completely reducible faithful $L$-module and $L<G$, so that by induction we conclude that $|L:L^c|<|V|\ (**)$.
Let $d\in L-K$ and $L_i=K_i\langle d \rangle$. Clearly $d$ acts on each of the $K_i$, and hence we see that
$L^c=[L, L, L]\leq \prod_{i=1}^{p} [L_i, L_i, L_i]\ =\ \prod_{i=1}^{p} L_i^c$.
As $L_i^c\leq K_i\cap D^c\ =\ K_i^c$ (by $(*)$) and obviously $K_i^c\leq L_i^c$, altogether we have $L_i^c=K_i^c$ for all $i$,
and thus $L^c=\prod_{i=1}^{p} K_i^c=K^c$. Hence with $(**)$ we get that
\[|V|>|L:L^c|=|L:K^c|=p|K:K^c|>p\cdot \frac{|V|}{p}=|V|,\]
which is a contradiction. This proves that indeed $D=K$.\\

We finally have accumulated sufficiently detailed information on the structure of $G$ so that now we can
finish the proof quite quickly.\\

Recall that $G^c=D^c=K^c$. Let $g\in G-D$. Let $k_1\in K_1-K_1^c$. Clearly we may assume
that $V_1^g=V_2$. Then $K_1^g=K_2$ and so $k_1^g\in K_2-K_2^c$. Thus
\[ [k_1, g]=k_1^{-1}k_1^g\in (K_1\times K_2)-(K_1^c\times K_2^c),\]
so \[1\ne [k_1, g]K^c=[k_1K^c,gK^c]\in (G/G^c)'\leq Z(G/G^c), \] since $G/G^c$ is of class 2.
Since
\[ (k_1^{-1}K_1^c)(k_1^gK_2^c)\prod\limits_{i=3}^{p}K_i^c\ =\ [k_1, g]K^c\in Z(G/G^c)\cap D/D^c\leq Z(D/D^c)=\X_{i=1}^{p}Z(K_i/K_i^c),\]
we conclude that $k_1^{-1}K_1^c\in Z(K_1/K_1^c)$ and thus
$k_1K_1^c\in Z(K_1/K_1^c)$, and as $k_1\in K_1-K_1^c$ was
chosen arbitrarily, we obtain that $K_1/K_1^c$ is abelian. Thus $K_1'\leq K_1^c\leq K_1'$ and so $K_1'=K_1^c$
which implies that $K_1$ is abelian. Hence all the $K_i$ and therefore $K$ is abelian,
and so $1=K^c=G^c$. (From $[k_1, g]K^c\in Z(G/G^c)$ it also follows immediately that $p=2$, but we do not need this any more.)
Since $G^c=1$, $G$ is now of class 2, and we obtain the desired conclusion $|G|<|V|$
from Glauberman's result in \cite[Proposition 1]{glauberman},
and this completes the proof of the theorem.
\end{beweis}

We finally can prove the inequality $|G/G^c| <|F(G)/F(G)^c|$ as announced above.\\

\begin{ncor}\label{cor1}
Let $G$ be a finite solvable group and $V$ a finite faithful completely reducible $G$-module,
possibly of mixed characteristic. If $G$ is not nilpotent, then
\[|G/G^c|\ <\ |F(G)/F(G)^c|.\]
\end{ncor}
\begin{beweis}
Suppose that $G$ is not nilpotent and proceed exactly as in the
beginning of the proof of \ref{thm1} where it is proved that $G$ can be assumed to be nilpotent. Note that in that proof the
inequality (1) there is obtained by using an inductive hypothesis which we do not have if we just want to prove the
current corollary. But we can, instead, simply use \ref{thm1} to get that inequality (1) and otherwise just follow the argument
outlined in the proof of \ref{thm1} to arrive at the inequality (4) which is the desired one.
\end{beweis}

Remark: We construct a few examples to show that our result is the best possible.

1. $G=D_8$ acts on $V={F}_3^2$, $D_8$ is of class $2$ and the bound is tight.

2. $G=SD_{16}$ (the Sylow $2$-subgroup of $\GL(2,3))$ acts on $V={F}_3^2$, and $SD_{16}$ is of class $3$ and has order $16$.

3. $G=Z_3 \wr Z_3$ acts on $V={F}_2^6$, $Z_3 \wr Z_3$ is of class $3$ and has order $81$.\\

\centerline{\bf Acknowledgements \rm}

The first author is partially supported by a grant from the Simons Foundation (\#280770 to Thomas M. Keller). The second author is partially supported by a grant from the Simons Foundation (\#499532 to Yong Yang).



\begin{thebibliography}{99}
\bibitem{aschbacher} M. Aschbacher, On the maximal subgroups of the finite classical groups, Invent. Math. {\bf 76}
(1984), 469--514.
\bibitem{aschbacher-guralnick} M. Aschbacher, R. M. Guralnick, On abelian quotients of primitive groups, Proc. Amer. Math. Soc. {\bf 107}
(1989), 89--95.
\bibitem{glauberman} G. Glauberman, On Burnside's other $p^{a}q^{b}$ theorem, Pacific J. Math. {\bf 56} (1975), 469--476.
\bibitem{hargraves}  B. B. Hargraves, The existence of regular orbits for nilpotent groups, J. Algebra {\bf 72} (1981), 54--100.
\bibitem{huppert} B. Huppert, Endliche Gruppen I, Springer, Berlin, 1967.
\bibitem{keller-yang} T. M. Keller, Y. Yang, Abelian quotients and orbit sizes of solvable linear groups,  Israel J. Math. {\bf 211} (2016), 23--44.
\bibitem{manz-wolf} O. Manz, T. R. Wolf, Representations of solvable
groups, London Mathematical Society, Lecture Notes Series {\bf 185}, Cambridge University Press,
Cambridge, 1993.
\bibitem{robinson-thompson} G. R. Robinson, J. G. Thompson, On Brauer's $k(B)$-problem, J. Algebra {\bf 184} (1996), 1143--1160.
(2002), 95--113.



\end{thebibliography}
\end{document}